\documentclass[11pt, reqno]{amsart}
\usepackage{amssymb,verbatim,amscd,amsmath,amssymb, bm}
\usepackage{mathrsfs}
\usepackage{braket}
\usepackage{array}
\usepackage{mathtools}
\usepackage{verbatim}
\usepackage{caption}

\usepackage{subcaption}
\usepackage{pdfsync}
\usepackage{tabu}
\usepackage{fullpage}
\usepackage{color}
\usepackage[svgnames,table]{xcolor}
\usepackage{enumerate}
\usepackage{pst-node,pst-text,pst-3d,pstricks}
\usepackage[dvips,all]{xy}
\usepackage[pagebackref]{hyperref}
\usepackage{graphics}
\usepackage{tikz}
\usepackage{pst-plot}
\usepackage{caption,subcaption}
\usetikzlibrary{patterns}
\usetikzlibrary{calc}
\usepackage{cleveref}
\usepackage{nicematrix}

\usetikzlibrary{matrix}
\usetikzlibrary{arrows}
\usetikzlibrary{decorations.pathmorphing}
\usepackage{graphicx,supertabular,paralist, colortbl}
\usetikzlibrary{decorations.pathreplacing, decorations.markings}
\tikzset{->-/.style={decoration={markings,mark=at position #1 with {\arrow{>}}},postaction={decorate}}}

\usepackage{tikz-cd}
\numberwithin{equation}{section}
\newtheorem{theorem}{Theorem}[section]
\newtheorem{lemma}[theorem]{Lemma}
\newtheorem{proposition}[theorem]{Proposition}
\newtheorem{corollary}[theorem]{Corollary}

\theoremstyle{definition}
\newtheorem{definition}[theorem]{Definition}

\newtheorem{def-prop}[theorem]{Definition-Proposition}
\newtheorem{remark}[theorem]{Remark}
\newtheorem{example}[theorem]{Example}

\newtheorem*{acknowledgment*}{Acknowledgment}
\newtheorem{setting}[theorem]{Setting}
\newtheorem{notation}[theorem]{Notation}

\newtheorem{algorithm}[theorem]{Algorithm}

\DeclareMathOperator{\Null}{Null}

\DeclareMathOperator{\gens}{gens}

\DeclareMathOperator{\supp}{supp}

\definecolor{navyblue}{rgb}{0.0, 0.0, 0.5}
\definecolor{darkred}{rgb}{0.55, 0.0, 0.0}

\newcommand{\D}{\mathcal{D}}
\newcommand{\C}{\mathcal{C}}

\def\P{{\mathcal P}}

\def\a{{\bf a}}

\def\u{{\bf u}}

\def\w{{\bf w}}
\def\x{{\bf x}}

\def\1{{\bf 1}}
\def\0{{\bf 0}}

\newcommand*\fixitem {\item[]%
  \refstepcounter{enumi}\hskip-\leftmargin\labelenumi\hskip\labelsep}

\begin{document}

\title{Toric Ideals of Weighted Oriented Graphs}

\author{Jennifer Biermann}
\address{Department of Mathematics and Computer Science, Hobart and William Smith Colleges \\
300 Pulteney St.
Geneva, NY 14456, USA}
\email{biermann@hws.edu}
\urladdr{}

\author{Selvi Kara}
\address{Department of Mathematics, University of Utah, 155 1400 E, Salt Lake City, UT 84112, USA}
\email{selvi@math.utah.edu}
\urladdr{}

\author{Kuei-Nuan Lin}
\address{Department of Mathematics,
Penn State University Greater Allegheny\\ 4000 University Dr, McKeesport, PA 15132, USA}
\email{kul20@psu.edu}
\urladdr{}

\author{Augustine O'Keefe}
\address{Department of Mathematics and Statistics, Connecticut College\\
270 Mohegan Avenue Pkwy,
New London, CT 06320, USA}
\email{aokeefe@conncoll.edu}
\urladdr{}

\begin{abstract} 
Given a vertex-weighted oriented graph, we can associate to it a set of monomials. We consider the toric ideal whose defining map is given by these monomials. We find a generating set for the toric ideal for certain classes of graphs which depends on the combinatorial structure and weights of the graph. We provide a result analogous to the unweighted, unoriented graph case, to show that when the associated simple graph has only trivial even closed walks, the toric ideal is the zero ideal. Moreover, we give necessary and sufficient conditions for the toric ideal of a weighted oriented graph to be generated by a single binomial and we describe the binomial in terms of the structure of the graph. 
\end{abstract}

\maketitle


\section{Introduction}

Toric ideals are of general interest due to a multitude of applications to other fields. One such application is using toric methods for analyzing chemical reaction networks, \cite{CDSS, Ga}. Given a connected chemical reaction network one can associate to it a monomial ideal which, in turn, gives rise to a toric ideal called the complex balancing ideal of the toric dynamical system. See \cite{CLS} for a more recent treatment from the commutative algebraic viewpoint of multi-Rees algebras.

Within the field of algebraic statistics, toric ideals arise as the vanishing ideal of discrete exponential families, also known as toric models. Generating sets of the toric ideals give information about Markov chains on sets of contingency tables \cite{DS}. As such, a generating set of the toric ideal is called a Markov basis. It is also known that the degree of the toric ideal gives rise to an upper bound for the maximum likelihood degree of the associated toric model \cite{CHKS}.

Toric ideals are defined via the kernel of a monomial map. When the monomials defining the map are square-free, one can use the methods of toric ideals of (hyper)graphs \cite{OH,PS,V}. In general, however, the defining monomials of the toric ideals in these applications are not square-free. Recently, it has been shown that the combinatorial structure of general toric ideals can be encoded by those arising from square-free monomials \cite{PTV}. However, the fundamental problem of finding an explicit generating set for toric ideals defined by non-square-free monomials is still open.

In this paper, we consider toric ideals arising from the edge ideals of vertex-weighted oriented graphs (VWOGs) which themselves have applications in coding theory, specifically in Reed-Muller-type codes \cite{MPV}. The study of ideals associated to VWOGs is very recent, and consequently, we have few results in this direction. In the case of monomial edge ideals of VWOGs, Cohen-Macaulayness has been studied \cite{GMSV,HLMRV,pitones2019monomial}, and the authors calculated the regularity of classes of edge ideals of VWOGs \cite{BBLO}.

The goal of this work is to characterize a generating set for the ideal based on the combinatorial structure of the graph. In the case of unoriented simple graphs, a generating set of the associated toric ideal is given by the even closed walks in the graph \cite{OH, V}. In the case of VWOGs, not every even closed walk gives an element of the toric ideal. Additionally, we give examples where the toric ideal of a VWOG has minimal generators corresponding to odd walks in the graph which would not happen in the unoriented case (see \Cref{ex:twoCycle} and \Cref{5 mingen}). So we see that the situation for the toric ideals of VWOGs is much more complicated than that of the unoriented case. 

The remainder of this paper is structured as follows. We first prove a result, \Cref{ZERO}, that is analogous to the unoriented, unweighted case that shows that if the radical of a monomial ideal is associated to a simple graph that has only trivial even closed walks then the toric ideal is the zero ideal. The simplest non-zero toric ideals of weighted oriented graphs are therefore even cycles.  We give a combinatorial description (which we call \emph{balanced}, see \Cref{balanced}) of when a weighted oriented even cycle has a non-zero toric ideal in (\Cref{balancedCycle}). Unlike the edge ideal case, a weighted oriented graph that consists of two cycles connected via a vertex or a path always has a non-zero toric ideal regardless of the parity of the lengths of the cycles. In \Cref{OneGen}, building off of the results on VWOGs with two cycles in \Cref{twoCycles}, we characterize when the toric ideal of a VWOG has a single generator.

Our main results are summarized in the following.

\begin{theorem}
\begin{enumerate}
   \fixitem Let $M$ be a monomial ideal such that $\sqrt{M}$ is generated by square-free monomials of degree 2 and its associated graph has only trivial even closed walks, then the toric ideal associated to $M$ is a zero ideal. 
 \item Let $\mathcal{D}$ be a weighted oriented cycle on $n$ vertices. The toric ideal associated to $\mathcal{D}$, $I_{\mathcal{D}}$, is non-zero if and only if  $\mathcal{D}$ is balanced.
  \item Let $\D$ be a weighted oriented graph. We provide necessary and sufficient conditions for which associated toric ideal $I_{\D}$ is generated by a single element. 
\end{enumerate}
\end{theorem}

\section{Preliminaries}\label{Preliminaries}
\subsection{Toric Ideals and Rings}

Let $M$ be a monomial ideal in polynomial ring $k[x_1,\dots,x_{\nu}]$ with minimal generating set $\gens(M)=\{\mathbf{x}^{\a_1},\mathbf{x}^{\a_2},\dots, \mathbf{x}^{\a_\mu}\}$, and $k[e_1,\dots,e_\mu]$ a polynomial ring with as many indeterminates as elements of $\gens(M)$. One can then define a $k$-algebra homomorphism $\varphi:k[e_1,\dots,e_{\mu}]\rightarrow k[x_1,\dots,x_{\nu}]$ via the monomial map $e_i\mapsto \mathbf{x}^{\a_i}$. The kernel of $\varphi$ is the \textbf{toric ideal of $M$}, denoted $I_M$, and the \textbf{toric ring of $M$} is the image of $\varphi$, denoted $k[M]$. Note that $k[M]$ has a natural multigrading by $\mathbb{N}\{\a_1,\a_2,\dots,\a_\mu\}$. Let $\a_j=(a_{1,j}, \dots, a_{\nu,j}) \in \mathbb{N}^\nu$, and let $A(M)=(a_{i,j})$ be the $\nu \times \mu$ incidence matrix of $M$ with $\mathscr{A}=\{\a_1,\dots,\a_\mu\}$ the set of columns in $A(M)$.   It is well known that $I_\mathcal{D}$ is generated by irreducible binomials $\prod_{i=1}^\mu e_i^{r_i}-\prod_{i=1}^\mu e_i^{s_i}$ where the monomials $\prod_{i=1}^\mu e_i^{r_i}$ and $\prod_{i=1}^\mu e_i^{s_i}$ have the same multidegree in $\mathbb{N}\mathscr{A}$, i.e. $\sum_{i=1}^\mu r_i\a_i=\sum_{i=1}^\mu s_i\a_i$ and $r_i$ and $s_i$ cannot both be non-zero for any $i$. For a particular element $f = \prod_{i=1}^\mu e_i^{r_i}-\prod_{i=1}^\mu e_i^{s_i}$ we say the {\bf support} of $f$, denoted $\supp f$, is the set of variables $e_i$ which occur in $f$ with non-zero exponent.

\subsection{Edge Ideals of Graphs and Weighted Oriented Graphs}
Let $G=(V(G),E(G))$ be a simple graph with vertex set $V(G)=\{x_1,\dots,x_\nu\}$ and edge set $E(G)=\{e_1=\{x_{1_1},x_{1_2}\},\dots,e_\mu=\{x_{\mu_1},x_{\mu_2}\}\}$. We adopt the definitions of cycle, tree, induced subgraph, leaf, walk, closed walk, and trivial closed walk from the literature of simple graphs. We recall some of the definitions here to be used in this work's content. 

\begin{definition}
A \emph{walk of length $q$} of a graph $G$ connecting $x_{i_1}\in V(G)$ and $x_{i_{q+1}}\in V(G)$ is a finite sequence of the form \[ W=(\{x_{i_1},x_{i_2}\},\{x_{i_2},x_{i_3}\},\dots,\{x_{i_q},x_{i_{q+1}}\}) \] with each $\{x_{i_k},x_{i_{k+1}}\}\in E(G)$. An \emph{even (resp. odd) walk} is a walk of even (resp. odd) length. A walk $W$ is called \emph{closed} if $x_{i_1}=x_{i_{q+1}}$, i.e. the initial vertex is the same as the ending vertex. We say $W$ is a \emph{path} if  $x_{i_j}\neq x_{i_k}$ for all $1\leq j<k \leq q$. A  \emph{cycle} is a closed path. A \emph{trivial closed walk} is a closed walk having no induced cycles.  
\end{definition}

Let $M=(x_{i_1}x_{i_2} |\{x_{i_1},x_{i_2}\}\in E(G))$ be a square-free monomial ideal associated to the graph $G$. We call $M$, the \emph{edge ideal of} $G$. Conversely, given a monomial ideal, $M$, such that $\gens (M)$ is a set of degree 2 square-free monomials, $x_{1_1}x_{1_2},\dots,x_{{\mu}_1}x_{\mu_2}\in k[x_1,\dots,x_{\nu}]$, then one can construct a simple graph $G_M$ with vertex set $\{x_1,\dots,x_\nu\}$ and edge set $\{\{x_{1_1},x_{1_2}\},\dots,\{x_{\mu_1},x_{\mu_2}\}\}$. The following is the well-known characterization of the generators of $I_M$ found in \cite{OH} and \cite{V}.

\begin{theorem}\label{GraphEvenClosedWalk}
When $M$ is the square-free monomial ideal associated to a simple graph $G$, the toric ideal $I_M$ is generated by even closed walks of $G$. 
\end{theorem}

\begin{example}
Let $\gens(M)=\{x_1x_2,x_2x_3,x_3x_1,x_1x_4,x_4x_5,x_5x_1\}$, then $G_M$ is comprised of two 3-cycles joined at a vertex with $E(G)=\{e_1=\{x_1,x_2\},e_2=\{x_2,x_3\},\dots,e_6=\{x_5x_1\}\}$. The shortest nontrivial even closed walk in $G_M$ traverses each edge of the graph once, and $I_M=(e_1e_3e_5-e_2e_4e_6)$.
\end{example}

\begin{remark}
\Cref{GraphEvenClosedWalk} shows that if the graph of a square-free monomial ideal, $M$, has only trivial even closed walks, then the toric ideal $I_M$ is the zero ideal.
\end{remark}

Let $\mathcal{D}=(V(\mathcal{D}), E(\mathcal{D}))$ be a weighted oriented graph with weight vector $\w\in  \mathbb{N}^{|V(\D)|}$ where $V(\mathcal{D})=\{x_1,\dots,x_{\nu}\}$, $E(\mathcal{D})=\{e_1,\dots,e_{\mu}\}$, and $\w=(w_1,w_2,\dots,w_\nu)$. We abuse notation and define the polynomial rings $k[E(\mathcal{D})]=k[e_1,\dots,e_{\mu}]$ and $k[V(\mathcal{D})]=k[x_1,\dots,x_{\nu}]$. We then consider the $k$-algebra homomorphism $\varphi: k[e_1,\dots,e_{\mu}]\rightarrow k[x_1,\dots,x_{\nu}]$ defined by $\varphi(e_i)=x_{i_1}x_{i_2}^{w_{i_2}}$ 
where $e_i$ is the directed edge from $x_{i_1}$ to $x_{i_2}$.

\begin{notation}
    We use the ordered pair $(x_{i_1},x_{i_2})$ to denote the oriented edge from $x_{i_1}$ to $x_{i_2}$. We write $e=\{x_{i_1}, x_{i_2}\}$ when the oriented edge $e$ is incident to $x_{i_1}$ and $x_{i_2}$ but the orientation is unspecified. 
\end{notation}
 
Let $A(\mathcal{D})=(a_{i,j})$ be the $\nu\times \mu$ incidence matrix of $\mathcal{D}$ defined by 
\[
    a_{i,j}=
    \begin{cases}
        1 & \text{ if $e_j=(x_i,x_k)\in E(\mathcal{D})$ for some $1\leq k\leq \nu$},\\
        w_i & \text{ if $e_j=(x_k,x_i)\in E(\mathcal{D})$ for some $1\leq k\leq \nu$},\\
        0 & \text{ otherwise}
    \end{cases}
\]
with $\mathscr{A}=\{\a_1,\dots,\a_{\mu}\}$ the set of columns in $A(\mathcal{D})$.
\\


Note that the monomials defining the $k$-algebra homomorphism $\varphi$ are exactly the generators of the monomial edge ideal studied in \cite{BBLO, GMSV, HLMRV, MPV, pitones2019monomial}. Defining the toric edge ring of $\mathcal{D}$ as the image of the monomial map defined by the generators of the monomial edge ideal is a natural analogue to defining the toric edge ring of an unoriented finite simple graph (see \cite{OH, V}).

\begin{remark}\label{rem:disjoint}
If $\mathcal{D}=\mathcal{D}_{1}\cup\mathcal{D}_{2}$ such that $V(\mathcal{D}_{1})\cap V(\mathcal{D}_{2})=\emptyset$,
then $I_{\mathcal{D}}=I_{\mathcal{D}_1}+ I_{\mathcal{D}_2}$. Thus we consider graphs with only one connected component.
\end{remark}

 In a manner similar to that used in the case of unoriented graphs, we characterize generators of the toric ideal of weighted oriented graphs in terms of cycles and closed walks. Our characterizations depend on the orientation of the edges in the oriented graph, and the cycle structure of its underlying unoriented graph. As such, when we refer to a cycle in a weighted oriented graph, we are considering a cycle in the underlying unoriented graph. When we also want to consider the orientation of the edges on the cycle, we use the term oriented cycle. An oriented cycle is not the same as a directed cycle (a term commonly used in the theory of directed graphs), where the latter requires the edges to all point in the same direction.
\section{Toric zero ideals}

In this section, we extend \Cref{GraphEvenClosedWalk}, namely if a simple graph has no non-trivial even closed walk then the toric ideal is the zero ideal.  We start with the main setting of this section which is more general than in the remainder of the paper.  Here we are associating any monomial $\x^{\a_i}=x_{i_1}^{a_{i_1,i}}x_{i_2}^{a_{i_2,i}}$ with $i_1\neq i_2$  to the edge $\{x_{i_1},x_{i_2}\}$ of a graph.

\begin{setting}\label{set:graph}
Let $M\subseteq k[x_1,\dots,x_\nu]$ be a monomial ideal with minimal generating set $\gens(M)=\{\bf x^{\bf a_1},\dots,\bf x^{\bf a_\mu}\}$ such that $|\supp {\bf x}^{{\bf a}_i}|=2$ for all ${\bf x}^{{\bf a}_i}\in \gens(M)$, and whenever ${\bf x}^{{\bf a}_i},{\bf x}^{{\bf a}_j}\in \gens(M)$ such that $i\neq j$ then $\supp{({\bf x}^{{\bf a}_i})}\neq \supp{({\bf x}^{{\bf a}_j})}$. Note that $\sqrt{M}$ is then a square-free monomial ideal of degree 2. We define the \emph{graph of $M$}, denoted $G(M)$, to be the graph whose edge ideal is $\sqrt{M}$. We take $E(G(M))=\{e_1,\dots,e_\mu\}$ so that the $I_M$ is the toric ideal defined by $M$ in $k[e_1,\dots,e_\mu]$.
\end{setting}

\begin{proposition}\label{lem:trimTrees}
Let $M$ and $G(M)$ be defined as in \Cref{set:graph}. Furthermore suppose that $M=M_1+M_2$ where $|V(G(M_1))\cap V(G(M_2))|=1$ and $G(M_2)$ is a tree. Then $I_M=I_{M_1}$. 
\end{proposition}

\begin{proof}
Let $f= \prod_{i=1}^{\mu}e_{i}^{r_{i}}-\prod_{i=1}^{\mu}e_{i}^{s_{i}}$ be an irreducible element in $I_M$. Let $E$ be the collection of all edges in the support of $f$ that are in $G(M_2),$ and suppose $E\neq\emptyset$. Since $G(M_2)$ is a tree, the subgraph of $G(M_2)$ with the edge set $E$ must have a leaf, say $x_t.$ Let $e_\ell$ be the unique edge incident to $x_t$ in $E$.   Therefore, since  $\sum_{i=1}^\mu r_i\mathbf{a}_i=\sum_{i=1}^\mu s_i\mathbf{a}_i$, where $\mathbf{a}_1\dots, \mathbf{a}_\mu$ are the exponent vectors of the generators $M$, we must have $r_\ell=s_\ell\neq 0$.  This contradicts the irreducibility of $f$, and thus $E$ must in fact be empty.  In other words, the support of every irreducible generator of $I_{M}$ consists of edges of $G(M_1)$.  \\
In the case that $G(M)$ is a tree we can take $M_1$ to be generated by a single generator of $M$, thus giving us $I_{M}=I_{M_1}=\{0\}$.
\end{proof}

\begin{corollary}\label{cor:removeWhisker}
Again let $M$ and $G(M)$ be as in \Cref{set:graph}. Assume further that $M=M'+\sum_{i=1}^kM_i$ such that each $G(M_i)$ is a tree,  $|V(G(M'))\cap V(G(M_i))|=1$ for all $i=1,2,\dots,k$, and $V(G(M_i))\cap V(G(M_j))=\emptyset$ for all $i\neq j$. Then $I_M=I_{M'}$.
\end{corollary}

\begin{remark}\label{rem:noWhisker}
Since leaves do not contribute to the generators of the toric ideal by \Cref{cor:removeWhisker}, for the remainder of the paper, we may assume that $G(M)$ has no leaves. 
\end{remark}

\begin{proposition}\label{lem:oddcycle}
Let $M$ and $G(M)$ be as in \Cref{set:graph}.  Suppose that $G(M)$ is an odd cycle, then $I_M$ is the zero ideal.
\end{proposition}

\begin{proof}
 Suppose $f=\prod_{i=1}^{\mu}e_i^{r_i}-\prod_{i=1}^{\mu}e_i^{s_i}\neq0$ is an irreducible element of $I_M$.  For convenience, we set $f_+=\prod_{i=1}^{\mu}e_i^{r_i}$ and $f_-=\prod_{i=1}^{\mu}e_i^{s_i}$ and note that $\varphi(f_+)=\varphi(f_-)$. Let $H$ be the subgraph of $G(M)$ induced by the edges in the support of $f$. Note that $H$ must then either be a forest or $H=G(M).$ In the first situation, we must have $f=0$ by \Cref{lem:trimTrees}. Thus, we may assume that $H=G(M)$.  Since $f$ is irreducible we have $\supp (f_+)\cap \supp (f_{-})=\emptyset$. Since $H=G(M)$, we must have $\supp \varphi(f)=\{x_1,\dots,x_{2n+1}\}$ which in turn implies $\supp\varphi(f_+)=\supp\varphi(f_-)=\{x_1,\dots,x_{2n+1}\}$. In the cycle, each vertex $x_i$ is incident to exactly two edges, $e_{i-1}$ and $e_i$ (where we take $e_0=e_{2n+1}$) and so one of $e_{i-1}$ or $e_i$ is in $\supp(f_+)$ while the other in $\supp(f_-)$. Without loss of generality, we can take $e_1\in\supp(f_+)$ so that $e_2\in\supp(f_-)$ which in turn tells us $e_3\in\supp(f_+)$. Traveling around the cycle in this manner will lead us to $e_{2n+1}\in\supp(f_+)$. Since both edges incident to $x_1$ are then in $\supp(f_+)$, $x_1\not\in\supp\varphi(f_-)$ thus leading to a contradiction.
\end{proof}

We are now ready to state the main theorem of this section.
\begin{theorem}\label{ZERO}
Let $M$ and $G(M)$ be defined as in \Cref{set:graph} such that $G(M)$ has only trivial even closed walks. Then the toric ideal $I_M$ is the zero ideal.
\end{theorem}
\begin{proof}
Saying that $G(M)$ has only trivial even closed walks means that each connected component of $G(M)$ can contain at most one cycle, which must be of odd length. By \Cref{rem:disjoint}  and \Cref{cor:removeWhisker} we may then assume that $G(M)$ is an oriented odd cycle. Now the theorem follows from \Cref{lem:oddcycle}.
\end{proof}

\section{Unicyclic graphs}

In the remainder of this paper, we turn our attention to finding a generating set for the toric ideal of a weighted oriented graph as defined in \Cref{Preliminaries}. We start by analyzing those graphs with a single cycle.  This is not a straightforward generalization of the unoriented case as the orientations of the edges of the cycle and the weights affect the generators of the ideal. 

Let $\mathcal{C}_{n}$ be a weighted oriented graph whose underlying unoriented graph is a cycle of
length $n$. Then the incidence matrix $A(\C_n)$ is an $n\times n$
matrix of the following form: 
\[
A(\mathcal{C}_{n})=\left[\begin{array}{cccccccc}
a_{1,1} & 0 & & \dots & &  & 0 & a_{1,n}\\
a_{2,1} & a_{2,2}& & \dots & & & 0 & 0\\
0 & a_{3,2} & \ddots\\
 &  & \ddots & a_{i,i}\\
\vdots & \vdots  &  & a_{i+1,i}&&& \vdots &\vdots\\
 &  &  &  && \ddots\\
 0 & 0 &\dots  &  &  & \ddots& a_{n-1,n-1}&0\\
 0& 0 & \dots & & & 0 & a_{n,n-1} & a_{n,n}
\end{array}\right]
\]
where $a_{i,i-1}, a_{i,i}\in\{1,w_{i}\}$ under the convention that $a_{1,0} = a_{1,n}$. We observe that if $\det(A(\mathcal{C}_{n}))\neq0$,
then the matrix has full rank. This implies the toric ideal, $I_{\mathcal{C}_{n}}$
is zero. On the other hand, if $\det(A(\mathcal{C}_{n}))=\prod_{i=1}^{n}a_{i,i}+(-1)^{n+1} a_{1,n}\prod_{i=1}^{n-1}a_{i+1,i}=0$,
then $\text{Null}(A(\mathcal{C}_{n}))$ is non-zero and there is at least
one generator in $I_{\mathcal{C}_{n}}.$ Due to the cycle structure,
any generator in $I_{\mathcal{C}_{n}}$ must involve all edges of
the cycle which we now prove. 

\begin{lemma}\label{cycleAllsupport}
If $f$ is a non-zero element of $I_{\mathcal{C}_{n}}$, then $\supp (f)=E(\C_n)$.
\end{lemma}

\begin{proof}
Let $0\neq f\in I_{\C_n}$ and $\textbf{u}=(u_1,\ldots, u_{n})$ its corresponding non-trivial element in $\text{Null} (A(\C_n))$ so that $A(\C_n)\u=\sum_{i=1}^{n} u_i \textbf{a}_i= \textbf{0}$ 
 where $\textbf{a}_i$ is the $i^{\text{th}}$ column of $A(\mathcal{C}_n)$. We claim that $u_i \neq 0$ for all $i\in[n]$ from which the lemma immediately follows. 
 For the sake of contradiction suppose there exists $e_j \in  E(\C_n)$ such that $u_j=0$, and let $x_\ell$ be a vertex incident to edge $e_j$ in $\mathcal{C}_n$.  Up to a relabeling of the graph, we may assume that $\ell \neq n$. Then let $e_{j'}$ be the only other edge incident to $x_\ell$ in $\mathcal{C}_n$.  Since the $\ell^\text{th}$ component of the vector $A(\C_n)\u$ satisfies $(A(\C_n)\u)_\ell= \sum_{i = 1}^n u_i a_{\ell,i} = 0$, and by assumption $u_i = 0$ for all $i \neq j'$, we must have $u_{j'} = 0$ as well.  An iterative application of this argument results in $u_i=0$ for all $e_i \in  E(\mathcal{C}_n),$ which contradicts the fact that $\textbf{u}=(u_1,\ldots, u_{n})$ is a non-trivial element in $\text{Null} (A(\mathcal{C}_n))$.
\end{proof}

By \Cref{cycleAllsupport}, if $\det(A(\mathcal{C}_{n}))=0$, there exists a minimal generator $f$ of $I_{\mathcal{C}_{n}}$ such that $\supp f$ involves all edges of the cycle. 
In order to have $\det(A(\mathcal{C}_{n}))=\prod_{i=1}^{n}a_{i,i}+(-1)^{n+1} a_{1,n}\prod_{i=1}^{n-1}a_{i+1,i}=0$,
$n$ must be an even number and $\prod_{i=1}^{n}a_{i,i}=a_{1,n}\prod_{i=1}^{n-1}a_{i+1,i}$.
This motivates the following definition.

\begin{definition}\label{balanced}
If $\C_n$ is a weighted oriented cycle on $n$ vertices and $A(\C_n)$ is its incidence matrix, we say $\C_n$ is \emph{balanced} if $n$ is even and $\prod_{i=1}^{n}a_{i,i}=a_{1,n}\prod_{i=1}^{n-1}a_{i+1,i}$.
\end{definition}

\begin{theorem}\label{balancedCycle}
Let $\C_n$ be a weighted oriented cycle on $n$ vertices.  Then the toric ideal $I_{\C_n}$ is non-zero if and only if $\C_n$ is balanced.

\end{theorem}

\begin{proof}

As noted above, the toric ideal $I_{\C_n} \neq 0$ if and only if the incidence matrix $A(\C_n)$ does not have full rank, or equivalently if and only if $\det(A(\C_n)) = 0$.  Therefore $I_{\C_n} \neq 0$ if and only if 
\[
\prod_{i = 1}^n a_{i,i} +(-1)^{n+1}a_{1,n}\prod_{i = 1}^{n-1} a_{i+1, i} = 0.
\]
Or equivalently
\[
\prod_{i = 1}^n a_{i,i} =(-1)^{n}a_{1,n}\prod_{i = 1}^{n-1} a_{i+1, i} .
\]
Since all $a_{i,i}, a_{i+1,i}, a_{1,n} \geq 1$, this last equation can hold if and only if $n$ is even and $\prod_{i=1}^{n}a_{i,i}=a_{1,n}\prod_{i=1}^{n-1}a_{i+1,i}$.
\end{proof}

The next result follows immediately from either \Cref{balancedCycle} or \Cref{lem:oddcycle}.

\begin{corollary}\label{oddCycle}
If $\C_n$ is a weighted oriented cycle on an odd number of vertices, then $I_{\C_n} = 0$. 
\end{corollary}

The following algorithm produces the single generator of the toric ideal of a balanced cycle.

\begin{algorithm}\label{removejoint2}
Input: $\C_n$ a balanced oriented weighted cycle with edges $e_1, \dots, e_n$ with weight vector $\textbf{w} = (w_1\dots w_n)$.  The algorithm outputs the generator of $I_{\C_n}$.  

\begin{itemize}

\item [Step 0:] Set $i=1$ and $r_1 = 1$.

\item [Step 1:] For $2 \leq i \leq n$, 

\indent if $e_{i-1} = (x_i, x_{i-1})$ and $e_i = (x_{i+1}, x_i)$, then $r_{i} = \frac{r_{i-1}}{w_{i}}$

\indent else if $e_{i-1} = (x_{i-1}, x_i)$ and $e_i = (x_i, x_{i+1})$, then $r_{i} = w_{i}r_{i-1}$

\indent else $r_i=r_{i-1}$.

\item [Step 2:] If no $r_i$ is fractional, then set $b_i=r_i$ for $1\leq i\leq n$

\indent else, $b_i=d*r_i$ where $d$ is the least common multiple of the denominators of $r_1,...,r_n$ in reduced form.
 
\end{itemize}

Output: $e_1^{b_1}e_3^{b_3}\dots e_{n-1}^{b_{n-1}} - e_2^{b_2}e_4^{b_4}\dots e_n^{b_n}$.  

\end{algorithm}

\begin{remark}

Note that the output binomial of \Cref{removejoint2} is indeed an element of the toric ideal. The construction forces the exponent vector of each monomial of $\varphi(e_1^{b_1}e_3^{b_3}\dots e_{n-1}^{b_{n-1}} - e_2^{b_2}e_4^{b_4}\dots e_n^{b_n})$ to have the same value at each variable. 
The positions of the non-zero entries of the incidence matrix $A(\C_n)$ guarantee that the rank of  $A(\C_n)$ is at least $n-1$, whereas the definition of a balanced cycle requires $\det(A(\C_n))=0$.  Hence the nullspace of $A(\C_n)$ has dimension 1.  Since we obtained the vector $(b_1, \dots, b_n)$ by multiplying by the least common multiple of the denominators, every integer-valued vector in the nullspace is a scalar multiple of $(b_1, \dots, b_n)$.   
\end{remark}

\begin{example}  \label{8cycles}

Let $\D_1$ be the weighted 8-cycle with edges $(x_1, x_2), (x_2, x_3), (x_3, x_4), (x_5, x_4), (x_6, x_5),\\ (x_7, x_6), (x_8, x_7), (x_1, x_8)$ and weight vector, $\w=(1,4,4,3,2,2,2,2)$.
The toric ideal of $\D_1$ is
\[
I_{\D_1} = (e_1e_3^{16}e_5^8e_7^2 - e_2^4e_4^{16}e_6^4e_8).
\]

Let $\D_2$ be the weighted 8-cycle with edges $(x_1, x_2), (x_2, x_3), (x_3, x_4), (x_5, x_4), (x_6, x_5),\\ (x_7, x_6), (x_8, x_7), (x_1, x_8)$ and weight vector, $\w=(1,8,2,3,2,2,2,2)$. 
The toric ideal of $\D_2$ is
\[
I_{\D_2} = (e_1e_3^{16}e_5^8e_7^2 - e_2^8e_4^{16}e_6^4e_8).
\]

Let $\D_3$ be the weighted 8-cycle with edges $(x_1, x_2), (x_2, x_3), (x_4, x_3), (x_5, x_4), (x_6, x_5), \\(x_7, x_6), (x_8, x_7), (x_1, x_8)$ and weight vector, $\w=(1,48,4,3,2,2,2,2)$. 

The toric ideal of $\D_3$ is
\[
I_{\D_3} = (e_1e_3^{48}e_5^8e_7^2 - e_2^{48}e_4^{16}e_6^4e_8).
\]

The weighted oriented graphs $\D_1$, $\D_2$, and $\D_3$ are all balanced. The graphs $\D_1$ and $\D_2$ differ only by the weights on $x_2$ and $x_3$. We see this difference manifests in their respective generators, specifically the exponent on $e_2$. The differences between $\D_1$ and $\D_3$ are the orientation of the edge on $\{x_3,x_4\}$ and the weight of $x_2$. Here we see this difference in the generators via the exponents on $e_3$ and $e_2$. 
\end{example}

From  \Cref{8cycles}, we see that in general, the condition for a cycle to be balanced is a complex interaction between the weights on the vertices and the orientations of the edges. We consider the extreme cases (when the edges are oriented in the same direction and when the weights on the vertices are all the same) in the next two corollaries.

We first consider the case where all of the edges of the cycle are oriented in the same direction. Formally, we say that a cycle 
$\{x_{1},x_{2}\}, \{x_{2},x_{3}\}, \dots \{x_{n},x_{1}\}$ is \emph{naturally oriented} if the orientations on the edges are either all $(x_i,x_{i+1})$ or all $(x_{i+1},x_i)$, where we take $x_{n+1}=x_1$. 

\begin{corollary}\label{lem:naturalCycle}
Let $\D$ be a weighted oriented graph with a unique cycle that is naturally oriented on $n$ vertices with $w_i \geq 2$ for some $i\in[n]$. Then $I_{\D}$ is the zero ideal.
\end{corollary}
\begin{proof}
First, suppose that $\D$ consists of a single cycle.  If $\mathcal{D}$ is an odd cycle, then by \Cref{oddCycle} $I_{\mathcal{D}}$ must be the zero ideal regardless of the orientation of the edges.  Let $\mathcal{D}$ be a naturally oriented even cycle. Then $\det A(\D)= \prod_{i=1}^{n}a_{i,i}-a_{1,n}\prod_{i=1}^{n-1}a_{i+1,i} = \pm(1 - \prod_{i=1}^{n}w_i)\neq 0,$ as $w_i> 1$ for some $i\in[n]$. Thus $\D$ is not balanced, and by \Cref{balancedCycle} $I_{\mathcal{D}}$ must be the zero ideal.  The full result then follows from \Cref{cor:removeWhisker}. 
\end{proof}

We next consider oriented cycles for which the weights are all the same.  For this, we will want to consider graphs in which half the edges are oriented in a clockwise manner and half counter-clockwise.  To make this precise we introduce the following definition.

\begin{definition}\label{def:balancedCycle}
Let $\mathcal{C}_n$ be a weighted oriented graph with its underlying graph the cycle on $n$ vertices.  Assume that $n$ is even. Let $E(\mathcal{C}_n)=\{e_1,\dots,e_{n}\}$ such that for all $i=1,\dots,n$, $e_i=\{x_i,x_{i+1}\}$ where we take $x_{n+1}=x_1$.
We then define the following subsets of $E(\mathcal{C}_n)$.\\
\begin{align*}
    \mathcal{C}_n^{inc}&=\{e_i=(x_i,x_{i+1})~:~i=1,\dots,n\}\\
    \mathcal{C}_n^{dec}&=\{e_i=(x_{i+1},x_i)~:~i=1,\dots,n\}=E(\mathcal{C}_n)\setminus\mathcal{C}_n^{inc}\\
\end{align*}
The oriented even cycle $\mathcal{C}_n$ is called \emph{uniformly balanced} if and only if $|\mathcal{C}_n^{inc}|=|\mathcal{C}_n^{dec}|$ and $w_i=w$ for all  $1\leq i \leq n$ such that the vertex $x_i\in V(\mathcal{C}_n)$ is neither a sink nor a source.
\end{definition}

\begin{corollary}\label{lem:even}
Let $\C_n$ be a weighted oriented cycle on $n$ vertices such that $w_i=w\geq 2$ for all non-source and non-sink $x_i\in V(\mathcal{C}_n)$. The toric ideal $I_{\C_n}$ is non-zero if and only if $\C_n$ is uniformly balanced.
\end{corollary}

\begin{proof}

Let $\C_n$ be a weighted oriented cycle on $n$ vertices in which $w_i = w \geq 2$ for all non-source and non-sink vertices $x_i$.  Let $A(\C_n)$ be the incidence matrix of $\C_n$.  By \Cref{balancedCycle}, $I_{\C_n} \neq 0$ if and only if $\prod_{i=1}^{n}a_{i,i}=a_{1,n}\prod_{i=1}^{n-1}a_{i+1,i}$ and $n$ is even. By assumption, $\C_n$ is uniformly balanced and thus $n$ is even, so we just need to verify the equality.  Notice that if $x_p$ is a sink vertex with weight $w_p$, the $p$-th row of $A
(\C_n)$ has $a_{p,p-1}=w_p$ and $a_{p,p}=w_p$. Hence the sink vertex does not impact the equality $\prod_{i=1}^{n}a_{i,i}=a_{1,n}\prod_{i=1}^{n-1}a_{i+1,i}$. Similarly, if $x_p$ is a source then $a_{p,p-1}=1$ and $a_{p,p}=1$.  By assumption all non-source and non-sink vertices have the same weight.  In the $i$-th column of the $A(\C_n)$ we have two possibilities:  either $a_{i,i} = w$ and $a_{i+1,i} = 1$ or vice versa, again taking $a_{n+1,n}=a_{1,n}$.  In the first case, the  edge $e_i=(x_{i+1},x_i) \in \C_n^{dec}$ 
and in the second $e_i=(x_i,x_{i+1}) \in \C_n^{inc}$.  Thus $\prod_{i=1}^{n}a_{i,i}=a_{1,n}\prod_{i=1}^{n-1}a_{i+1,i}$ if and only if half of the edges are in $\C_n^{inc}$ and half are in $\C_n^{dec}$.  Thus $I_{\C_n} \neq 0$ if and only if $\C_n$ is uniformly balanced. 
\end{proof}

\begin{example}  \label{8cycles2}

Let $\D_1$ be the weighted 8-cycle with edges $(x_1, x_2), (x_2, x_3), (x_4, x_3), (x_5, x_4), (x_6, x_5),\\ (x_6, x_7), (x_7, x_8), (x_1, x_8)$ and with all vertices with weight 3.  The toric ideal of $\D_1$ is then
\[
I_{\D_1} = (e_1^3e_3^9e_5e_7^3-e_2^9e_4^3e_6e_8^3).
\]
Let $\D_2$ be the weighted 8-cycle with edges $(x_1, x_2), (x_2, x_3), (x_4, x_3), (x_5, x_4), (x_5, x_6),\\ (x_6, x_7), (x_7, x_8), (x_1, x_8)$ also with all vertices having weight 3.  The toric ideal of $\D_2$ is then the zero ideal.
The only difference between $\D_1$ and $\D_2$ is the orientation of the edge on vertices $\{x_5,x_6\}$ so that $|\D_2^{inc}|\neq |\D_2^{dec}|$.

\end{example}

\section{Two Connected Cycles}\label{twoCycles}

We now consider generating sets of the toric ideals of connected oriented graphs comprised of distinct oriented cycles $\C_m$ and $\C_n$ with varying degrees of overlap. More specifically, we determine for which of these graphs the toric ideal can be generated by a single element. We then give a combinatorial description of the support of that generator.

Before we state the main theorem of this section, we introduce one additional notation. Let $\P_s$ denote an oriented path of length $s$, i.e. its underlying unoriented graph is a path. Note that the incidence matrix of $\P_s$ is of size $(s+1)\times s$ and has the following form: 
\[
A(\mathcal{P}_{s})=\left[\begin{array}{ccccccc}
a_{1,1} & 0 & & \dots & &  & 0 \\
a_{2,1} & a_{2,2}& & \dots & & & 0 \\
0 & a_{3,2} & \ddots\\
 &  & \ddots & a_{i,i}\\
\vdots & \vdots  &  & a_{i+1,i}&&& \vdots \\
 &  &  &  && \ddots\\
 0 & 0 &\dots  &  &  & \ddots& a_{s,s}\\
 0& 0 & \dots & & & 0 & a_{s+1,s} 
\end{array}\right].
\]

\begin{theorem}\label{thm:singleGen} 
    Let $\D$ be a weighted oriented graph comprised of two oriented cycles $\C_m$ and $\C_n$ satisfying one of the following:
    \begin{enumerate}[(A)]
        \item $\C_m$ and $\C_n$ share a single vertex,
        \item the set $E(\C_m)\cap E(\C_n)$ induces an oriented path with at least one edge, or
        \item $\C_m$ and $\C_n$ are connected by an oriented path $\P_s$ of length $s\geq 1$. 
    \end{enumerate}
    Then the toric ideal $I_\D$ is generated by a single element if and only if at most one of its cycles is balanced. 
\end{theorem}

\begin{figure}[h]
    \begin{subfigure}{0.48\textwidth}
	    \centering
	    \includegraphics[scale=0.9]{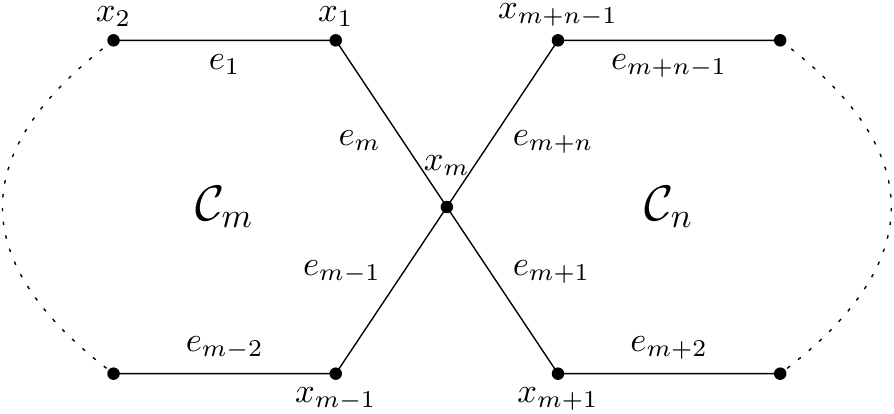}
	    \caption{Two cycles sharing a vertex}
	    \label{fig:twoCyclesGV}
	\end{subfigure}
	\begin{subfigure}{0.48\textwidth}
        \centering
	     \includegraphics[scale=0.9]{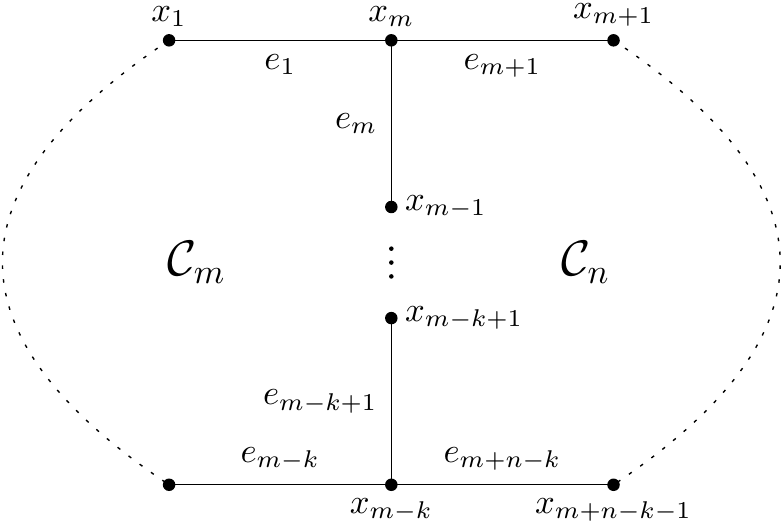}
	     \caption{Two cycles sharing edges}
	    \label{fig:twoCyclesGE}
	\end{subfigure}
	
	\bigskip
	
	\begin{subfigure}{0.8\textwidth}
		\centering
	    \includegraphics[scale=0.9]{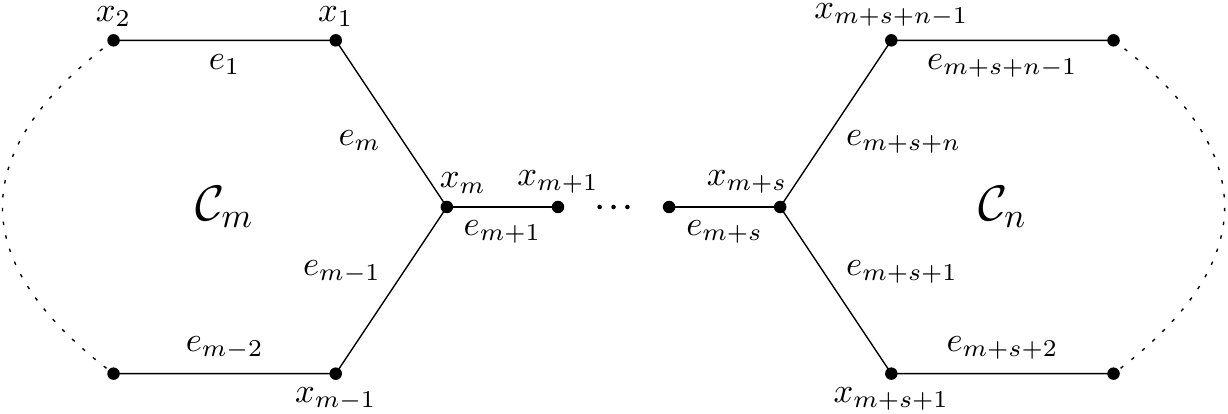} 
	    \caption{Two cycles connected by a path}
	    \label{twoCyclesGP}
    \end{subfigure}
    \caption{Cycle forms from \Cref{thm:singleGen}}
    \label{fig:cycleForms}
\end{figure}

 \Cref{fig:cycleForms} illustrates the three possible forms of the underlying unoriented cycle structure of $\D$ along with a labeling of vertices and edges that will be used throughout this section.

    \begin{figure}[h]
    \centering
    \begin{subfigure}{0.48\textwidth}
	    \centering
	    \includegraphics[scale=1]{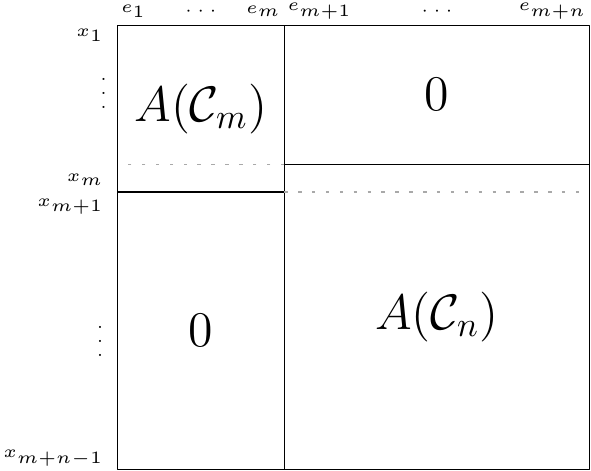}
	    \caption{Two cycles sharing a vertex}
	    \label{fig:twoCyclesV}
	\end{subfigure}
    \begin{subfigure}{0.48\textwidth}
		\centering
	    \includegraphics[scale=1]{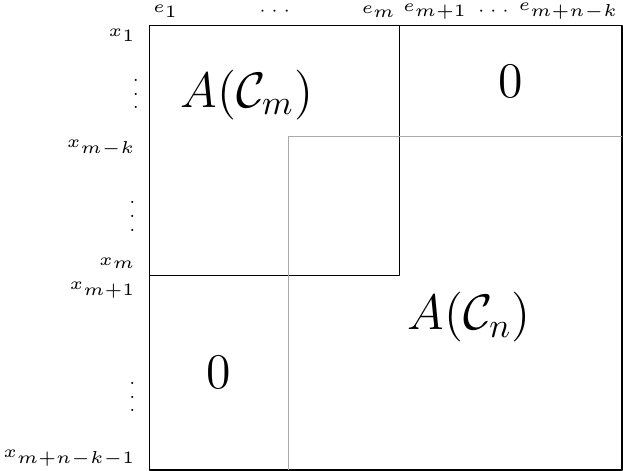} 
	    \caption{Two cycles sharing edges}
	    \label{twoCyclesE}
    \end{subfigure}
	
	\bigskip
	
	\centering
	\begin{subfigure}{0.48\textwidth}
        \centering
	     \includegraphics[scale=1]{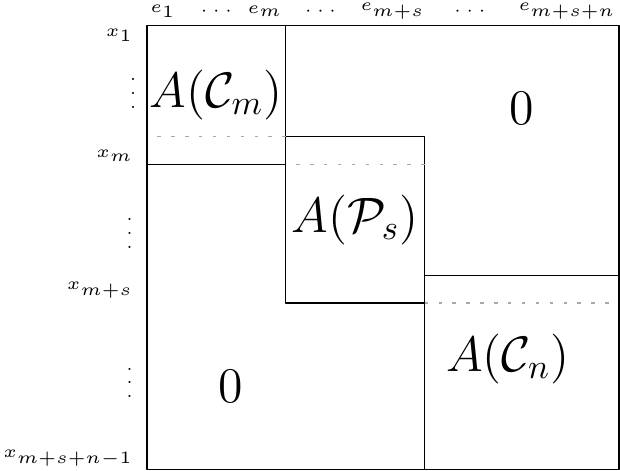}
	     \caption{Two cycles connected by a path}
	    \label{fig:twoCyclesP}
	\end{subfigure}
    \caption{Incidence matrix forms from \Cref{thm:singleGen}}
    \label{fig:matrixForms}
\end{figure}

 \begin{remark}\label{rem:threeCycles}
 Note that in \Cref{fig:cycleForms}(B) there are three distinct cycles to consider in $\D$: the two ``inner cycles" $\C_m$ and $\C_n$, and the ``outer" cycle obtained by deleting the edges $\{e_{m-k+1},\dots,e_m\}$. For $I_\D$ to be generated by a single element, at most one of these three cycles can be balanced. 
 \end{remark}

\begin{proof} 
    
    If both $\C_m$ and $\C_n$ are  balanced then \Cref{lem:even} tells us that their corresponding toric ideals are nonzero.  Using a straightforward linear algebra argument, one can extend a basis element of $\Null A(\C_m) $ and a basis element of $\Null A(\C_n)$ to linearly independent elements of $\Null A(\D)$. Thus $I_\D$ would have at least two generators. The same argument would apply in the case of overlapping edges for which two of the three cycles in $\D$ are balanced, see \Cref{rem:threeCycles}. 
    
    We now assume that $\D$ has at most one  balanced cycle, and furthermore suppose that $\C_m$ is not balanced. Throughout each of the three cases, we refer to the $(i,j)^\text{th}$ entry of $A(\D)$ as $a_{i,j}$. We also define $\widetilde{A(\D)}$ to be the square matrix obtained from the incidence matrix $A(\D)$ by deleting the last column. Informed by the labeling of vertices and edges as in \Cref{fig:cycleForms}, the matrices in \Cref{fig:matrixForms} illustrate the possible forms of the incidence matrix of $\D$. Note that in all three situations, $A(\D)$ has one more column than row so that $\widetilde{A(\D)}$ is square. Thus if $\det \widetilde{A(\D)}\neq 0$ then $A(\D)$ is full rank and $\dim_\mathbb{Q}\Null A(\D)=1$. This in turn shows $I_\D$ is principal.
 
    We first consider the case that $\C_m$ and $\C_n$ are joined by a single vertex. Define $B$ to be the matrix obtained from $A(\C_n)$ by deleting both the first row and the last column. Then $\det \widetilde{A(\D)}=\det A(\C_m)\det B$. Since $\C_m$ is not balanced we must have that $\det A(\C_m)\neq 0$. One can also calculate
        \[
            \det B = \prod_{i=m+1}^{m+n-1}a_{i,i}\neq 0
        \]
    where $a_{i,i}\in\{1,w_i\}$ for all $i$.
    Thus $A(\D)$ is of full rank and $I_\D$ has a single generator. 
    
    We now consider the case that $E(\C_m)\cap E(\C_n)$ induces a path of length $k\geq 1$. Note that in the block formation of $A(\D)$ as illustrated in \Cref{fig:matrixForms}(B) the entries below the overlap, i.e. the $(i,j)$-entries with $i\geq m+1$ and $m-k\leq j\leq m$, are all zero. Therefore, $\det \widetilde{A(\D)}=\det A(\C_m)\det B_k$
    where $B_k$ is obtained from $A(\C_n)$ by deleting its first $k$ columns, its first $k+1$ rows, and the last column. It can then be observed that 
    \[\det B_k =\prod_{i=m+1}^{m+n-k-1}a_{i,i}\neq 0.\]
     Again we see that $A(\D)$ is indeed of full rank. 
    
    Lastly, assume $\C_m$ and $\C_n$ are connected by an oriented path $\P_s$ with $s\geq 2$. Again define the matrix $B$ to be the matrix obtained from $A(\C_n)$ by deleting the first row and last column. We also take $E$ to be the matrix obtained from $A(\P_s)$ by deleting its first row. Note that the $(i,j)^\text{th}$ entry of $A(\P_s)=0$ if and only if $j>i+1$ or $j<i$. Thus $\det \widetilde{A(\D)}=\det A(\C_m)\det E\det B \neq 0$ since 
        \[
            \det E = \prod_{i=m+1}^{m+s}a_{i,i}\neq 0.
        \]
    Yet again we see that $A(\D)$ is of full rank. 
\end{proof}

In the following corollary, we give a combinatorial description of the generator guaranteed by \Cref{thm:singleGen} in terms of the edges appearing in the support of the generator. Again we will see that there are three cycles to consider in the case that $\C_m$ and $\C_n$ share edges. 

\begin{corollary}\label{cor:genForm}
     Let $\D$ be an oriented graph comprised of two cycles as in \Cref{thm:singleGen} such that $I_\D$ is generated by a single irreducible binomial $f$. Then $\supp f=E(\C)$ if one of the cycles $\C$ is balanced. Otherwise, $\supp f=E(\D)$. 
\end{corollary}

\begin{proof}
Let $\C_m$ and $\C_n$ be distinct cycles in $\D$ where $I_\D$ is generated by a single nonzero element $f$. By \Cref{thm:singleGen} we know at least one of these cycles is unbalanced, assume it is $\C_n$. If $\C_m$ is balanced, then we can trivially extend $\u\in\Null A(\C_m)$ to an element of $\Null A(\D)$ so that by \Cref{cycleAllsupport} we can assume the generator $f$ of $I_\D$ satisfies $\supp f=E(\C_m)$. Now assume that $\C_m$ is also unbalanced.

We first argue that $E(\C_m)\setminus E(\C_n)\subseteq\supp f.$ Let $\u\in\Null A(\D)$ so that $\sum_{i=1}^{|E(\D)|}u_i\a_i=\bf 0$, and for the sake of contradiction suppose $u_j=0$ for some $e_j\in E(\C_m)\setminus E(\C_n)$. Since $\D$ has no multiple edges, we can find $x_\ell\in V(\C_m)\setminus V(\C_n)$ incident to $e_j$. There is then only one other edge incident to $x_\ell$ which we can assume to be $e_{i+1}\in E(\C_m)$. It follows then that $u_{i+1}=0$. Iterating this argument yields $u_i=0$ for all $e_i\in E(\C_m)\setminus E(\C_n)$ which would imply that $\supp f\subseteq E(\D)\setminus E(\C_m)$. Note that the only cycle in $\D\setminus\C_m$ is $\C_n$ so that $I_\D=I_{\C_n}$ by \Cref{lem:trimTrees}. However, since $\C_n$ is unbalanced \Cref{balancedCycle} yields $I_\D=I_{\C_n}$ is the zero ideal, a contradiction. Since $\C_m$ is also unbalanced, a symmetric argument tells us we also have $E(\C_n)\setminus E(\C_m)\subseteq\supp f$.

In the case that $\C_m\cap\C_n$ is a single vertex we have 
    \[
        E(\C_m)\bigtriangleup E(\C_n)=E(\C_m)\cup E(\C_n)=E(\D)
    \]
where $E(\C_m)\bigtriangleup E(\C_n)$ denotes the symmetric difference of $E(\C_m)$ and $E(\C_n)$. Thus we must have $\supp f = E(\D)$. 
    
When $\C_m$ and $\C_n$ are connected by a path of length at least one, we again see that 
    \[
        E(\C_m)\bigtriangleup E(\C_m)=E(\C_m)\cup E(\C_n)\subseteq\supp f.
    \]
Let $\D$ be labeled as in \Cref{fig:cycleForms}(C). Note that since $\C_m$ and $\C_n$ are unbalanced, $\u\in\Null A(\D)$ would have to satisfy $u_i\neq 0$ for some $m+1\leq i\leq m+s$, otherwise $I_{\C_m}\neq (0) \neq I_{\C_n}$. In the case that $s=1$, then we would immediately have $\supp f= E(\D)$. Suppose that $s\geq 2$ and $u_j\neq 0$ for some $e_j\in\P_s$. Up to a relabeling of $\D$ we can assume that $x_j\in V(\P_s)\setminus[V(\C_m)\cup V(\C_n)]$ for which the only other edge incident to it is $e_{j+1}\in E(\P_s)$. Since $\sum u_i\a_i=\bf{0}$ and $u_j\neq 0$ we must then have $u_{j+1}\neq 0$. Iterating this argument results in $E(\P_s)\subseteq\supp f$ so that again $\supp f=E(\D)$.  

Lastly, we consider the case that $E(\C_m)\cap E(\C_n)$ induces a path of length $k$, call it $\P_k$. If the outer cycle induced by $E(\C_m)\bigtriangleup E(\C_n)$ is balanced, we can then assume that $f$ is supported by this cycle. If the outer cycle is unbalanced, then we must have $u_i\neq 0$ for some $e_i\in\P_k$. By a similar argument as above, we can then conclude that $u_i\neq 0$ for every $e_i\in\P_k$ which again results in $f$ being supported by all of the edges in $\D$. 

\end{proof}

\begin{remark}
    When constructing the toric ideal of an unoriented graph, one can construct a binomial generator $f$ from an even closed walk by assigning alternating edges along the walk to each monomial in $f$. One can similarly think of constructing the generator of the toric ideal $I_\D$ where $\D$ is of one of the forms in \Cref{thm:singleGen}.  However, the parity of the lengths of the cycles will affect in which of the monomials each edge appears. Furthermore, due to the weights on the vertices, it is difficult to calculate the exponents of each variable.
\end{remark}

The following example illustrates different cases of \Cref{thm:singleGen} and \Cref{cor:genForm}. We further discuss how the generating set of the toric ideal differs from that of the underlying graph.

\begin{figure}[h]
    \centering
    \begin{subfigure}{0.48\textwidth}   
        \centering
        \includegraphics[scale=0.9]{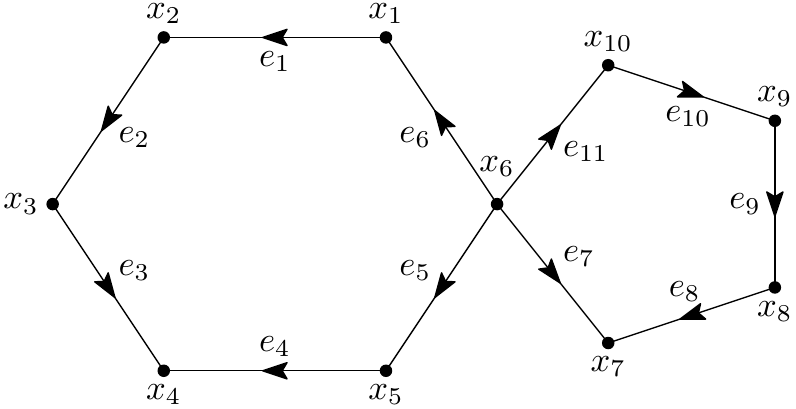}
        \caption{$\D_1$: Two unbalanced cycles}
    \end{subfigure}
    \begin{subfigure}{0.48\textwidth}
        \centering
        \includegraphics[scale=0.9]{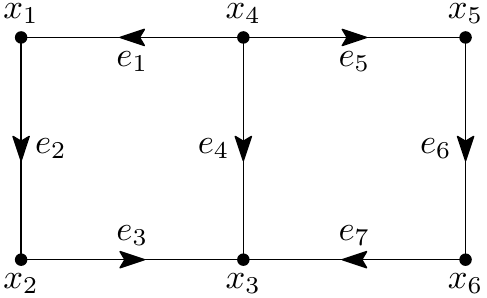}
        \caption{$\D_2$: One balanced cycle}
    \end{subfigure}
    
    \bigskip
    
    \begin{subfigure}{0.8\textwidth}
        \centering
        \includegraphics[scale=0.9]{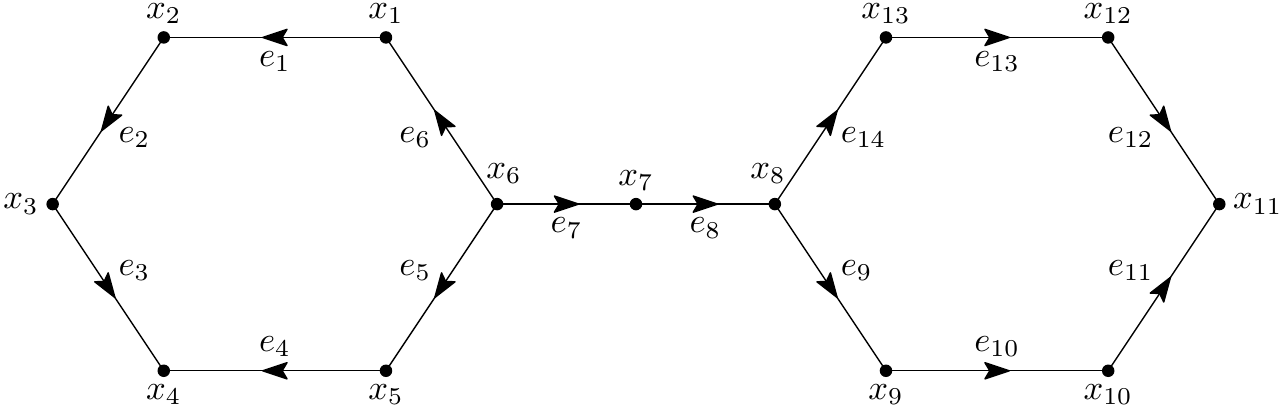}
        \caption{$\D_3$: One balanced cycle}
    \end{subfigure}
    \caption{Graphs with two cycles}
    \label{twoCycleGraph}
\end{figure}

\begin{example}\label{ex:twoCycle} 

We see in Figure \ref{twoCycleGraph}(A) $\D_1$ is comprised of an oriented $\C_6$ and an oriented $\C_5$ connected by a vertex. With the weight vector $\w=(2,2,3,4,6,1,4,2,3,2)$, neither of these cycles are balanced.  In this case
\[
    I_{\D_1}=(e_2^{52}e_4^{156}e_6^{13}e_7^{12}e_9^6e_{11}-e_1^{26}e_3^{156}e_5^{26}e_8^{12}e_{10}^2).
\]
Here the toric ideal of the underlying graph would also be generated by a single element, but it would come from the 6-cycle, not the whole graph.

In Figure \ref{twoCycleGraph}(B), $\D_2$ is an oriented $\C_4$ and an oriented $\C_4$ sharing one edge with weight vector $\w=(2,2,3,1,2,2)$. Both 4-cycles are unbalanced but the outer 6-cycle is balanced. Thus the 6-cycle admits the single generator,
\[
    I_{\D_2}=(e_1e_3^4e_6^2-e_2^2e_5e_7^4).
\]
In the case of the underlying graph, there would be two generators, each coming from one of the 4-cycles in the graph. Of course, the outer 6-cycle would still admit a binomial in the toric ideal, but it would not be needed as a minimal generator.  

In Figure \ref{twoCycleGraph}(C), $\D_3$ is comprised of two oriented 6-cycles connected by a path of length 2. With the weight vector $\w=(2,2,3,4,12,1,5,3,2,3,5,4,2)$ only the left cycle is balanced, and thus gives us the generator of $I_{\D_3}$. 

\[
    I_{\D_3}=(e_2^4e_4^{12}e_6-e_1^2e_3^{12}e_5)
\]
In contrast, the toric ideal of the underlying graph of $\D_3$ would have two minimal generators that could be taken to be from each of the two 6-cycles.

\end{example}

\section{Toric ideals generated by a single element}

The following theorem is stated as generally as possible, allowing for trees to be attached to the graph that can be removed without affecting the toric ideal by \Cref{lem:trimTrees}.  Note that the graph structures described in \Cref{OneGen} are the same as those in \Cref{thm:singleGen}.

\begin{theorem}\label{OneGen}
Let $\D$ be a weighted oriented graph for which the toric ideal $I_D$ is generated by a single element.  Then $\D$ is one of the following:
\begin{enumerate}
\item a balanced unicyclic graph,
\item a graph with exactly two oriented cycles connected by either a vertex or a path such that at most one of those cycles is balanced, or
\item a graph with exactly three cycles such that the intersection of any two of the cycles is a connected path and at most one of those cycles is balanced.
\end{enumerate}
\end{theorem}

\begin{proof}
Let $\D$ be a weighted oriented graph whose toric ideal $I_{\D}$ is generated by a single element.  Let $A(\D)$ be the incidence matrix of $\D$.  In order for the null space of $A(\D)$ to have dimension 1, we must have $|E(\D)| - |V(\D)| \leq 1$. 

Let $T$ be a spanning tree of $\D$.  It is well-known that $|E(T)|-|V(T)| = -1$. We know that $T \neq \D$ since the toric ideal of a weighted oriented tree is 0. Therefore the difference between $\D$ and $T$ is either one or two edges.  If $\D$ has one more edge than $T$, then $\D$ is a unicyclic graph.  

If $\D$ has two more edges than $T$, then as above adding the first edge to $T$ gives a unicyclic graph.  Let $\{v,w\}$ be the second edge.  If $v$ and $w$ are both on the cycle, then $\D$ is a graph with two cycles that share an edge.  If $v$ is on the cycle and $w$ is not then $\D$ is two cycles sharing some number of edges or two cycles sharing a vertex.  If $v$ and $w$ are both not on the cycle then $\D$ is either two cycles sharing some number of edges or two cycles sharing a vertex or two cycles connected by a bridge.  
\end{proof}

We conclude with an example that illustrates the complex behavior of the generators of toric ideals of more complicated weighted oriented graphs.

\begin{example}\label{5 mingen}
 
Let $\D_1$ and $\D_2$ be the weighted oriented graphs with 3 cycles connected by a vertex as shown in \Cref{fig:3cycles}. 

\begin{figure}[h]
    \centering
    
    \includegraphics[scale=0.85]{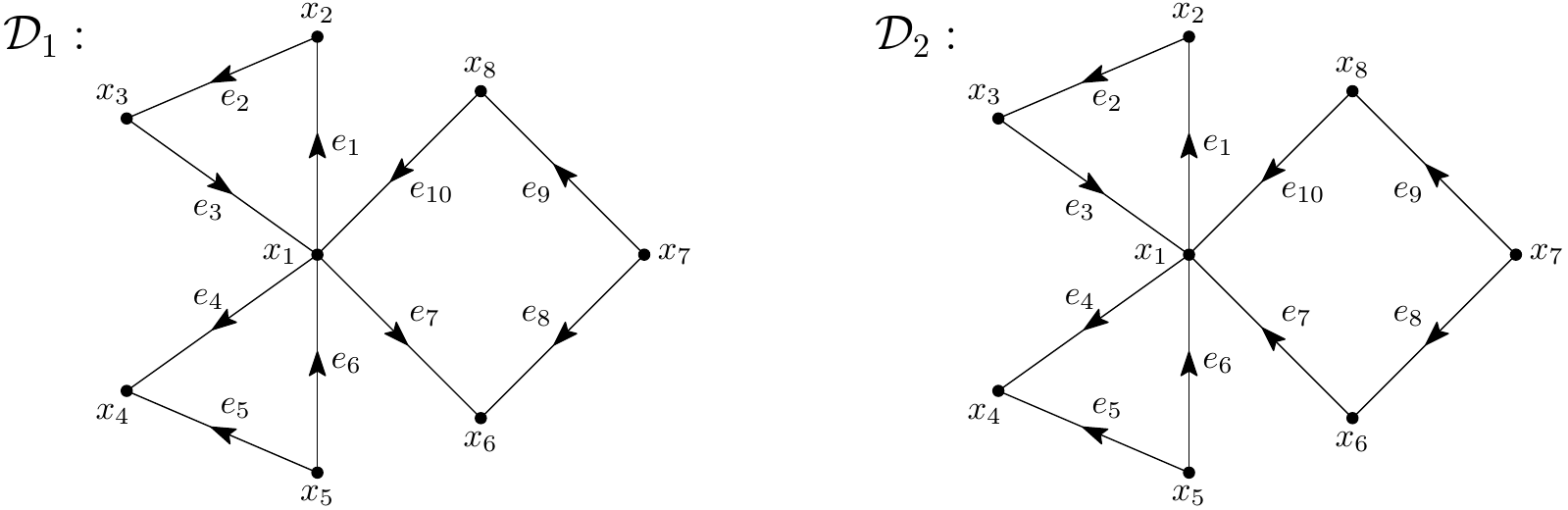}
    \caption{Three cycles joined at a vertex}
    \label{fig:3cycles}
\end{figure}

If $\D_1$ has weight vector $\w=(2,2,2,2,1,2,1,2)$, then the toric ideal of $\D_1$ is
\[
I_{\D_1} = \left(
    \renewcommand*{\arraystretch}{1.3}
    \begin{array}{lll}
    e_4e_6e_7e_9-e_5e_8e_{10}^2, & e_1e_3^4e_5^3-e_2^2e_4^3e_6^3,& e_1e_3^4e_5^2e_7e_9-e_2^2e_4^2e_6^2e_8e_{10}^2,\\ e_1e_3^4e_5e_7^2e_9^2-e_2^2e_4e_6e_8^2e_{10}^4, & e_1e_3^4e_7^3e_9^3-e_2^2e_8^3e_{10}^6
    \end{array}
    \right).
\]

When $\D_2$ has weight vector $\w=(2,2,2,2,1,2,1,2)$,
the toric ideal of $\D_2$ is
\[
I_{\D_2} = (e_7^2e_9-e_8e_{10}^2,e_1e_3^4e_5^3-e_2^2e_4^3e_6^3).
\]
The only difference between $\D_1$ and $\D_2$ is the orientation of the edge $\{x_1,x_6\}$. In $\D_2$ the 4-cycle is balanced but the 4-cycles is not in $\D_1$. We observe that any two cycles in $\D_1$ give a single generator, and there are two generators having all edges as their supports in $I_{\D_1}$. By contrast, $\D_2$ and the underlying graph of $\D_1$ only have two generators, one coming from the two 3-cycles, and the other from the 4-cycle. 

This example shows that even though the weights of the weighted oriented graphs are the same, the number of generators could vary wildly depending on the orientation of the edges. Furthermore, $\D_1$ has two generators having the same supports with different exponents which never happens in the simple edge graph case. 

\end{example}

\begin{acknowledgment*}
    The authors are grateful to the software system \texttt{Macaulay2} \cite{M2}, for serving as an excellent source of inspiration.
\end{acknowledgment*}

\bibliographystyle{plain}
\bibliography{references}
\end{document}